\documentclass[11pt,reqno]{amsart}

\theoremstyle{plain}
\newtheorem{theorem} {Theorem}[section]
\newtheorem{lemma}[theorem] {Lemma}
\newtheorem{proposition}[theorem] {Proposition}
\newtheorem{corollary}[theorem] {Corollary}

\theoremstyle{definition}
\newtheorem{definition}[theorem] {Definition}

\newtheorem{example} [theorem]{Example}

\theoremstyle{remark}
\newtheorem{remark}[theorem] {Remark}

\numberwithin{equation}{section}

\usepackage{amsfonts}
\usepackage{amssymb}
\usepackage{amscd}

\newcommand{\R}{{\mathbb R}}
\newcommand{\Z}{{\mathbb Z}}

\newcommand{\PP}{{\mathcal P}}

\newcommand{\C}{{\mathcal C}}
\newcommand{\CC}{{\mathbb C}}

\newcommand{\TT}{{\mathcal T}}

\newcommand{\B}{{\mathcal B}}
\newcommand{\A}{{\mathcal A}}
\newcommand{\BB}{{\mathfrak B}}
\newcommand{\al}{{\alpha}}
\newcommand{\la}{{\lambda}}
\newcommand{\sa}{{\sigma}}

\newcommand{\iy}{{\infty}}
\newcommand{\vphi}{{\varphi}}
\newcommand{\vep}{{\varepsilon}}
\newcommand{\g}{{\gamma}}

\newcommand{\bna}{\begin{eqnarray}}
\newcommand{\ena}{\end{eqnarray}}
\newcommand{\ba}{\begin{eqnarray*}}
\newcommand{\ea}{\end{eqnarray*}}
\newcommand{\beq}{\begin{equation}}
\newcommand{\eeq}{\end{equation}}

\begin{document}

\title[Bernstein--Markov Type Inequalities]
{A Note on Sharp Multivariate Bernstein-- and Markov--Type Inequalities}
\author{Michael I. Ganzburg}
 \address{212 Woodburn Drive, Hampton,
 VA 23664\\USA}
 \email{michael.ganzburg@gmail.com}

 \begin{abstract}
 Let $V$ be a symmetric convex body in $\R^m$.
 We prove sharp Bernstein--type inequalities for
entire functions of exponential type with the spectrum in $V$
and discuss certain properties of the extremal functions.
Markov--type inequalities with sharp constants
 for algebraic polynomials on $V$ and certain
 non-symmetric convex bodies  are proved as well.
 \end{abstract}
 \maketitle

 \section{Introduction}\label{S1}
\setcounter{equation}{0}
\noindent
 In this paper  we discuss several new
 Bernstein-- and Markov--type inequalities
 in the uniform norm
  with sharp constants for multivariate
 entire functions of exponential type and
 algebraic polynomials.
 In particular, we extend the Bernstein inequality
 to entire functions $f$ of exponential type
  with the spectrum
 in a symmetric convex body  $V\subset\R^m$
 in the form
 \ba
 \sup_{y\in K^*}\left\vert\sum_{j=1}^m\frac{\partial f(x)}
 {\partial x_j}y_j\right\vert\le M(K,V)\|f\|_{C(\R^m)},
 \qquad x\in\R^m,
 \ea
 where $K$ is
 a symmetric convex body in $\R^m$
 and $K^*$ is the polar of $K$.
 Formulae for the sharp constant $M(K,V)$
 in the above inequality and properties of its
 extremal functions are discussed as well.
 We also prove the similar Markov-type inequalities
 for polynomials on symmetric
 and certain non-symmetric convex bodies.
\vspace{.12in}\\
\textbf{Notation.}
Let $\R^m$ be the Euclidean $m$-dimensional space with elements
$x=(x_1,\ldots,x_m),\, y=(y_1,\ldots,y_m),
\,t=(t_1,\ldots,t_m),\,u=(u_1,\ldots,u_m)$,
the inner product $(t, y):=\sum_{j=1}^mt_jy_j$,
and the norm $\vert t\vert:=\sqrt{(t, t)}$.
Next, $\CC^m:=\R^m+i\R^m$ is the $m$-dimensional complex space with elements
$z=(z_1,\ldots, z_m)=x+iy$
and the norm $\vert z\vert:=\sqrt{\vert x\vert^2+\vert y\vert^2}$;
$\Z^m$ denotes the set of all integral lattice points in $\R^m$;
and $\Z^m_+$ is a subset of $\Z^m$
of all points with nonnegative coordinates.
We also use a multi-index $\al=(\al_1,\ldots,\al_m)\in \Z^m_+$
with
 $\vert\al\vert:=\sum_{j=1}^m\al_j$ and
 $x^\al:=x_1^{\al_1}\cdot\cdot\cdot x_m^{\al_m}$.

Let $C(E)$ be the space of all continuous complex-valued functions $F$
 on a measurable set $E\subseteq\R^m$  with the finite uniform norm
 $
 \|F\|_{C(E)}:=\sup_{x\in E} \vert F(x)\vert,
 $
 and let  $\PP_{n,m}$ be the set of all
  polynomials $P(x)=\sum_{\vert\al\vert\le n}c_\al x^\al$
   in $m$ variables of total degree at most $n,\,n\in Z^1_+$,
   with complex coefficients. In addition,
   $T_n(\cdot):=\cos(n\arccos \cdot)\in\PP_{n,1}$
 is the Chebyshev polynomial of the first kind.
 For a differentiable function $f:\R^m\to \CC^1$, let
 \ba
 \nabla f(x)
 :=\left(\frac{\partial f(x)}{\partial x_1},\cdots,
 \frac{\partial f(x)}{\partial x_m}\right),\qquad x\in\R^m,
 \ea
 be the gradient vector.

 Throughout the paper $R,\,R_1,\,R_2,\,C,\,C_0,\,C_1,\ldots$
 denote constants independent
of essential parameters.
Occasionally we indicate dependence on certain parameters.
 The same symbol does not
 necessarily denote the same constant in different occurrences.
 \vspace{.12in}\\
\textbf{Convex Bodies.} In this paper we need certain
 definitions and properties of
 convex bodies in $\R^m$,
 i.e., compact convex sets with non-empty interiors.
 We first define the \emph{width} $w(\C)$ of a
 convex body $\C$ in $\R^m$ as
 the minimum distance between two parallel supporting
  hyperplanes of $\C$
 and define the \emph{diameter} $d(\C)$ of $\C$ as
 the maximum distance between two points of $\C$.
 Let
$\BB^m:=\{t\in\R^m: \vert t\vert\le 1\},\,
Q^m:=\{t\in\R^m: \left\vert t_j\right\vert\le 1, 1\le j\le m\}$,
 and
$O^m:=\{t\in\R^m: \sum_{j=1}^m\left\vert t_j\right\vert\le 1\}$
be the $m$-dimensional ball, cube, and octahedron, respectively.

 Next, let $V$
  be a centrally symmetric (with respect to the origin)
 closed
 convex body in $\R^m$
 with the boundary $\partial(V)$, the width $w(V)$, and
 the diameter $d(V)$.
 In addition, let
 $V^*:=\{y\in\R^m: \forall\, t\in V, \vert (t, y)\vert \le 1\}$
 be the \emph{polar} of $V$.
 It is well known that $V^*$ is a centrally symmetric
  (with respect to the origin) closed
 convex body in $\R^m$ and $V^{**} =V$
 (see, e.g., \cite[Lemma 3.4.7]{SW1971}).
 The sets $V$ and $V^*$ generate the following norms
 on $\R^m$ and $\CC^m$ by
 \ba
  \|x\|_{V}
  :=\max_{t\in V^*}\vert (t, x)\vert,\,\,\,
  \|y\|_{V^*}
  :=\max_{t\in V}\vert (t, y)\vert,\,\,
  x, y\in\R^m;\quad
  \|z\|_{V^*}
  :=\max_{t\in V}\left\vert\sum_{j=1}^m t_jz_j\right\vert,\,\, z\in\CC^m.
 \ea
 Note also that $V$ and $V^*$ are the unit balls in the norms
 $\|\cdot\|_V$ and  $\|\cdot\|_{V^*}$ on $\R^m$, respectively.
 For example, the following convex bodies and their polars:
 \ba
 &&V_\mu:=\left\{x\in\R^m: \| x\|_{V_\mu}=
 \left(\sum_{j=1}^m\left\vert x_j\right\vert^{\mu}\right)^{1/\mu}\le 1\right\},\\
  &&(V_\mu)^*=\left\{y\in\R^m: \| y\|_{(V_\mu)^*}=
  \left(\sum_{j=1}^m\left\vert y_j\right\vert^\rho
 \right)^{1/\rho}\le 1\right\},
\ea
 where $\mu\in[1,\iy],\,\rho\in[1,\iy]$, and $1/\mu+1/\rho=1$,
 have various applications in analysis.
 In particular,
 $V_1=(V_\iy)^*=O^m,\,V_2=(V_2)^*=\BB^m,$ and $V_\iy=(V_1)^*=Q^m$.

 Throughout the paper $K$ and $V$ are
 centrally symmetric (with respect to the origin) closed convex bodies
 in $\R^m$. Next, we define the constant $M(K,V)$ and sets $\A$ and $\B$.
 \bna
 &&M(K,V):=\max_{y\in V}\|y\|_K
 =\max_{y\in \partial (V)}\|y\|_K
 =\max_{y\in\R^m\setminus \{0\}}\| y\|_K/\| y\|_V,\label{E1.1}\\
 &&\A=\A(K,V):=\{\mathbf{a}\in \partial (V):
 \|\mathbf{a}\|_K= M(K,V)\},\label{E1.1a}\\
 &&\B=\B(V^*,K^*):=\{\mathbf{b}\in \partial (K^*):
  \|\mathbf{b}\|_{V^*}= M(V^*,K^*)\}.\label{E1.1b}
 \ena
 Then the following property holds:
 \begin{proposition}\label{P1.1}
 For any $\mathbf{a}=
 \mathbf{a}(K,V)\in\A$ there exists $\mathbf{b}
 =\mathbf{b}\left(V^*,K^*\right)\in\B$ such that
 \beq\label{E1.2}
  M(K,V)=M\left(V^*,K^*\right)
  =\vert (\mathbf{b},\mathbf{a})\vert.
  \eeq
  \end{proposition}
  \proof
 Since
 \beq\label{E1.2b}
 \max_{y\in V}\| y\|_K
 =\max_{y\in V}\max_{x\in K^*}\vert(x,y)\vert
 =\max_{x\in K^*}\max_{y\in V}\vert(x,y)\vert
 =\max_{x\in K^*}\| x\|_{V^*},
 \eeq
 we arrive at the first equality in \eqref{E1.2}.
 Next, given $\mathbf{a}\in\A$, there exists
 $\mathbf{b}\in \partial\left(K^*\right)$ such that
 $ M(K,V)=\max_{y\in V}\|y\|_K=\|\mathbf{a}\|_K
 =\vert (\mathbf{b},\mathbf{a})\vert$.
 Since by \eqref{E1.2b},
 \ba
 \| \mathbf{b}\|_{V^*}=\max_{y\in V}\vert (\mathbf{b},y)\vert
 \le \max_{x\in K^*}\max_{y\in V}\vert(x,y)\vert
 =M(K,V)=\vert (\mathbf{b},\mathbf{a})\vert
 \le \| \mathbf{b}\|_{V^*},
 \ea
 we see that $\mathbf{b}\in\B$. Thus
 the second equality in \eqref{E1.2} is established.
 \hfill $\Box$

 Similarly one can prove that
 for any $\mathbf{b}\in\B$  there exists $\mathbf{a}\in\A$
 such that the second equality in \eqref{E1.2} is valid.
 More general forms of the first equality in \eqref{E1.2}
   are discussed in \cite[Theorem 3.13]{J2017}.
The following geometric characterization of $\mathbf{b}$
 in Proposition \ref{P1.1}
 follows from \eqref{E1.1a}, \eqref{E1.1b}, and \eqref{E1.2}:
 given $\mathbf{a}\in\A$, one can choose $\mathbf{b}\in\B$
 as a point,
 satisfying the condition
  $(1/M(K,V))\,\mathbf{b}
  \in \partial \left(V^*\right)\cap H_{\pm}(\mathbf{a})$,
  where $H_{\pm}(\mathbf{a})
 :=\{x\in\R^m:(\mathbf{a},x)=\pm 1\}$ are
 parallel supporting hyperplanes of $V^*$.
 Note that in certain cases
 $\mathbf{b}
=\pm\left(M\left(K,V\right)\vert
\mathbf{a}\vert^{-2}\right)
\mathbf{a}$
(see Example \ref{Ex2.5}).

  In particular, it follows from \eqref{E1.1}
  and \eqref{E1.2} that
  if $K$ is the ball $\BB^m$ with the norm
  $\|\cdot\|_K=\vert\cdot\vert$, then
  \beq\label{E1.2a}
  d\left(V^*\right)/2=M\left(\BB^m,V^*\right)
  =M\left(V,\BB^m\right)=2/w(V).
  \eeq
  \noindent
\textbf{Entire Functions of Exponential Type.}
 The set of all trigonometric polynomials
 $T(x)=\sum_{k\in V\cap \Z^m}c_ke^{i(k, x)}$ with complex
 coefficients is denoted by $\TT_{V}$.
 A more general set of entire functions is defined below.

  \begin{definition}\label{D1.1}
 We say that an entire function $f:\CC^m\to \CC^1$
  has exponential type $V$
 if for any $\vep>0$ there exists a constant $C_0(\vep,f)>0$ such that
 for all $z\in \CC$,
 $\vert f(z)\vert\le C_0(\vep,f)\exp\left[(1+\vep)\| z\|_{V^*}\right]$.
 \end{definition}
 \noindent
  The set of all entire function of exponential type $V$ is denoted
  by $B_V$.
  In case of $m=1$, we use the notation
  $B_\sa:=B_{[-\sa,\sa]},\,\sa>0$.
  Throughout the note, if no confusion may occur,
  the same notation is applied to
  $f\in B_V$ and its restriction to $\R^m$ (e.g., in the form
  $f\in  B_V\cap C(\R^m)$).
  The set $B_V$ was defined by Stein and Weiss
  \cite[Sect. 3.4]{SW1971}. For $V=\sa Q^m$ and
  $V=\sa \BB^m,\,\sa>0$, similar
  sets were
  defined by Bernstein \cite{B1948} and Nikolskii
  \cite[Sects. 3.1, 3.2.6]{N1969}, see also
  \cite[Definition 5.1]{DP2010}.

 The rest of the paper is organized as follows: In Section \ref{S2}
 we obtain Bernstein--type inequalities for functions from $B_V$,
 and in Section \ref{S4} we discuss
 Markov--type inequalities for polynomials  from $\PP_{n,m}$
 on symmetric and non-symmetric convex bodies.
 Certain properties of extremal functions in
 Bernstein--type inequalities are discussed in
 Section \ref{S3}.

 \section{Bernstein--type Inequalities}\label{S2}
\setcounter{equation}{0}
\noindent
Throughout the section a point
$\mathbf{a}\in \partial(V)$ satisfies the equality
$\left\| \mathbf{a}\right\|_K=M(K,V)$, i.e., $\mathbf{a}\in \A(K,V)$
(see \eqref{E1.1} and \eqref{E1.1a}),
and a point
$\mathbf{b}\in \B\left(V^*,K^*\right)$
 is defined by
\eqref{E1.1b} and \eqref{E1.2}.
\vspace{.12in}\\
\textbf{Bernstein Inequalities on $\R^1$.}
The Bernstein--Szeg\"{o} type inequality for
complex-valued functions
$\vphi\in B_\sa\cap C(\R^1), \sa>0$, can be presented in the following form:
\beq\label{E2.1}
\left\vert\sin \al\,\vphi^\prime(\tau)-\sa \cos \al\, \vphi(\tau)\right\vert
\le \sa\|\vphi\|_{C(\R^1)},\qquad \tau\in\R^1,\quad \al\in[0,2\pi).
\eeq
In particular, the classic Bernstein inequality
for
complex-valued functions
$\vphi\in B_\sa\cap C(\R^1), \sa>0$,
immediately follows
 from \eqref{E2.1} for $\al=\pi/2$,
\beq\label{E2.2}
\left\vert \vphi^\prime(\tau)\right\vert
\le \sa\|\vphi\|_{C(\R^1)},\qquad \tau\in\R^1.
\eeq
The functions
$
\vphi_0(\tau):=C_1 e^{i\tau}
+C_2 e^{-i\tau},\,C_1\in\CC^1,\, C_2\in\CC^1,
$
with $\|\vphi_0\|_{C(\R^1)}
=\left(\left\vert C_1\right\vert^2
+\left\vert C_1\right\vert^2\right)^{1/2}$
are the only \textit{extremal}
 functions in \eqref{E2.1} or \eqref{E2.2}, i.e.,
 equality holds in \eqref{E2.1} or \eqref{E2.2}
  for a certain
 $\tau=\tau_0\left(\sa,C_1,C_2\right)\in\R^1$
 and any $\al\in[0,2\pi)$ if and only if
 $\vphi=\vphi_0$.

For real-valued functions $\vphi\in B_\sa\cap C(\R^1)$,
\eqref{E2.1} can be reduced to the inequality
\beq\label{E2.3}
\left(\vphi^\prime(\tau)\right)^2+\sa^2 (\vphi(\tau))^2
\le \sa^2\|\vphi\|_{C(\R^1)}^2,\qquad \tau\in\R^1.
\eeq
The functions
$
\vphi_0(\tau):=R_1 \cos \sa\tau+R_2 \sin \sa\tau,
\,R_1\in\R^1,\, R_2\in\R^1,
$
with $\|\vphi_0\|_{C(\R^1)}
=\left(R_1^2+R_2^2\right)^{1/2}$
are the only extremal functions in \eqref{E2.3}, i.e.,
equality holds in \eqref{E2.3}
for a certain
 $\tau=\tau_0\left(\sa,R_1,R_2\right)\in\R^1$
if and only if
 $\vphi=\vphi_0$.

The proof of \eqref{E2.1} with the description of all extremal functions
in \eqref{E2.1} and \eqref{E2.3} can be found in \cite[Sect. 84]{A1965}
(see also \cite[Ch. 11]{B1954}).
\vspace{.12in}\\
\textbf{Bernstein--type Inequalities on $\R^m$.}
The following theorem extends inequality \eqref{E2.1} and
its extremal functions to multivariate entire functions
from $B_V\cap C(\R^m)$.

\begin{theorem}\label{T2.1}
For every $f\in B_V\cap C(\R^m),\,x\in\R^m,\,y\in
\partial (K^*)$, and
$\al\in[0,2\pi)$, the following inequality holds:
\beq\label{E2.4}
\left\vert\sin \al\,\left(\nabla f(x),y\right)
-M(K,V)\cos \al\, f(x)\right\vert
\le M(K,V)\|f\|_{C(\R^m)}.
\eeq
There exist $x=x_0\in\R^m$ and $y=y_0\in
\partial (K^*)$ such that
for any $\al\in[0,2\pi)$ and the functions
\beq\label{E2.5}
f(x)=f_0(x)
:=C_1 e^{i(\mathbf{a},x)}+C_2 e^{-i(\mathbf{a},x)},
\qquad C_1\in\CC^1,\quad C_2\in\CC^1,
\eeq
with $\left\|f_0\right\|_{C(\R^m)}
=\left(\left\vert C_1\right\vert^2
+\left\vert C_1\right\vert^2\right)^{1/2}$,
 equality holds in \eqref{E2.4}.
 Hence inequality \eqref{E2.4} is sharp.
\end{theorem}
\proof
\textbf{Step 1.}
We first need the following elementary lemma:
\begin{lemma}\label{L2.1a}
Given two complex numbers
$a+ib$ and $c+id$,
the function
\ba
h(\la):=\sup_{\al\in[0,2\pi)}
\left\vert\la\sin \al\,(a+ib)-\cos \al\,(c+id)\right\vert
\ea
is not decreasing on $(0,\iy)$.
\end{lemma}
\proof
If
$(a+ib)(c+id)= 0$, then the statement is trivial.
If $(a+ib)(c+id)\ne 0$, then
by a straightforward calculation for $\la\in(0,\iy)$,
\ba
&&h^2(\la)
=(1/2)\left(\la^2\left(a^2+b^2\right)+c^2+d^2
+\sqrt{\left[\la^2\left(a^2+b^2\right)-c^2-d^2\right]^2
+4\la^2(ac+bd)^2}\right),\\
&&h^\prime(\la)
= \frac{\la(a^2+b^2)}{2h(\la)}\left(1+
\frac{\la^2\left(a^2+b^2\right)-c^2-d^2
+4(ac+bd)^2}
{\sqrt{\left[\la^2\left(a^2+b^2\right)-c^2-d^2\right]^2
+4\la^2(ac+bd)^2}}\right)
\ge 0.
\ea
Note that the expression under the square roots above is
positive for $\la>0$. Thus the lemma is established.
\hfill $\Box$
\vspace{.12in}\\
\textbf{Step 2.}
Next, let $f\in B_V\cap C(\R^m),\,x\in\R^m$,
 and $y\in\R^m\setminus 0$,
and let
\ba
\vphi(w)=\vphi(\tau+i\g):=f(x+(\tau+i\g)y),
\qquad \tau\in\R^1,\quad \g\in\R^1,
\ea
be the restriction of $f(z)$ to the
 one-dimensional complex plane
$z=x+w y$ in $\CC^m$.
Then the following property holds:
\begin{lemma}\label{L2.1b}
$\vphi\in B_{\| y\|_{V^*}}\cap C(\R^1)$.
\end{lemma}
\proof
$\vphi$ is an entire function of $w\in\CC^1$, and for
any $\vep>0$,
\ba
\vert\vphi(w)\vert\le C_0(\vep,f)\exp\left[(1+\vep)\|
x+wy\|_{V^*}\right]
\le C_0(\vep,f)\exp\left[(1+\vep)\|x\|_{V^*}\right]
\exp\left[(1+\vep)\vert w\vert\,\| y\|_{V^*}\right].
\ea
Therefore, $\vphi$ is a univariate entire function of exponential type
$\| y\|_{V^*}$ and $\vphi\in C(\R^1)$.
\hfill $\Box$
\vspace{.12in}\\
\textbf{Step 3.}
Using now Step 2 and inequality \eqref{E2.1} for $\tau=0$, we obtain
\ba
\left\vert\sin \al\,\left(\nabla f(x),y\right)
- \| y\|_{V^*}\cos \al\, f(x)\right\vert
&=&\left\vert\sin \al\,\vphi^\prime(0)
- \| y\|_{V^*}\cos \al\, \vphi(0)\right\vert\\
&\le& \| y\|_{V^*}\,\|\vphi\|_{C(\R^1)}
\le \| y\|_{V^*}\,\|f\|_{C(\R^m)}.
\ea
Hence for any $x\in\R^m$ and $y\in \partial  (K^*)$,
\beq\label{E2.6}
\sup_{\al\in[0,2\pi)}\left\vert\sin \al\,\left(\nabla f(x),y\right)
\| y\|_{K^*}/\| y\|_{V^*}
- \cos \al\, f(x)\right\vert
\le \|f\|_{C(\R^m)}.
\eeq
Next, note that for any $y\in\R^m\setminus \{0\}$,
\beq\label{E2.7}
\| y\|_{K^*}/\| y\|_{V^*}\ge
 1/M(V^*,K^*)=1/M(K,V),
 \eeq
 by \eqref{E1.1} and \eqref{E1.2}.
 Finally applying Lemma \ref{L2.1a} for
 $a+ib=\left(\nabla f(x),y\right)$ and $c+id=f(x)$,
we obtain \eqref{E2.4} from \eqref{E2.6} and \eqref{E2.7}.
\vspace{.12in}\\
\textbf{Step 4.}
Finally, by Proposition \ref{P1.1},
for $\mathbf{a}\in \A$
there exists $\mathbf{b}=y_0\in \B$ such that
$(\mathbf{a},y_0)=\| \mathbf{a}\|_K=M(K,V)$.
Then the function $f_0$ defined by \eqref{E2.5}
belongs to $B_V\cap C(\R^m)$,
and there exists $x_0\in\R^m$ such that
\ba
&&\left\vert\sin \al\,\left(\nabla f_0(x_0),y_0\right)
-M(K,V)\cos \al f_0(x_0)\right\vert\\
&&=M(K,V)
\left\vert e^{-i(\al+\pi)}C_1 e^{i(\mathbf{a},x_0)}
+e^{i(\al+\pi)}C_2 e^{-i(\mathbf{a},x_0)}\right\vert
=M(K,V)\left\|f_0\right\|_{C(\R^m)}.
\ea
This completes the proof of Theorem \ref{T2.1}.
\hfill $\Box$\vspace{.12in}\\
Multivariate versions of \eqref{E2.2}
and \eqref{E2.3} are presented below.
\begin{corollary}\label{C2.2}
(a) For every $f\in B_V\cap C(\R^m)$ and $x\in\R^m$,
the following inequality holds:
\beq\label{E2.8}
\|\nabla f(x)\|_K
\le M(K,V)\|f\|_{C(\R^m)}.
\eeq
There exists $x=x_0\in\R^m$ such that
for functions \eqref{E2.5}
 equality holds in \eqref{E2.8}.
 Hence inequality \eqref{E2.8} is sharp.\\
 (b) For a real-valued function
 $f\in B_V\cap C(\R^m)$ and
 every  $x\in\R^m$,
 the following inequality holds:
\beq\label{E2.9}
\|\nabla f(x)\|_K^2+M^2(K,V) f^2(x)
\le M^2(K,V)\|f\|_{C(\R^m)}^2.
\eeq
There exists $x=x_0\in\R^m$ such that
for functions
\beq\label{E2.9a}
f(x)=f_0(x)
:=R_1 \cos \,(\mathbf{a},x)+R_2 \sin \,(\mathbf{a},x),
\qquad
R_1\in\R^1,\quad R_2\in\R^1,
\eeq
with $\left\|f_0\right\|_{C(\R^m)}
=\left(R_1^2
+R_2^2\right)^{1/2}$,
 equality holds in \eqref{E2.9}.
 Hence inequality \eqref{E2.9} is sharp.
\end{corollary}
\proof
(a) Choosing in \eqref{E2.4} $\al=\pi/2$
and $y\in \partial (K^*)$ such that
$\left\vert\left(\nabla f(x),y\right)\right\vert
=\|\nabla f(x)\|_K$,
we arrive at \eqref{E2.8}. In addition,
there exists $x_0\in\R^m$ such that the function $f_0$
defined by \eqref{E2.5} satisfies the equalities
\beq\label{E2.10}
\left\|\nabla f_0(x_0)\right\|_K
=\left\| \mathbf{a}\right\|_K
\left\vert C_1 e^{i(\mathbf{a},x_0)}
-C_2 e^{-i(\mathbf{a},x_0)}\right\vert
=M(K,V)\left\|f_0\right\|_{C(\R^m)}.
\eeq
(b) Inequalities  \eqref{E2.4} and
\beq\label{E2.11}
\left(\nabla f(x),y\right)^2
+M^2(K,V) f^2(x)
\le M^2(K,V)\|f\|_{C(\R^m)}^2,
\qquad x\in\R^m,\quad y\in \partial (K^*),
\eeq
 are equivalent for real-valued functions
 $f\in B_V\cap C(\R^m)$.
 Then choosing $y\in \partial (K^*)$ such that
$\left\vert\left(\nabla f(x),y\right)\right\vert
=\|\nabla f(x)\|_K$,
we obtain \eqref{E2.9} from \eqref{E2.11}.
 The case of equality in \eqref{E2.9} for functions $f_0$
  can be proved
 similarly to \eqref{E2.10}.
\hfill $\Box$

\begin{remark}\label{R2.3}
Inequalities \eqref{E2.1}, \eqref{E2.2},
 and \eqref{E2.3} are special
 cases of \eqref{E2.4}, \eqref{E2.8},
  and \eqref{E2.9}, respectively,
for $K=V=[-1,1]$.
The author \cite[Theorem 1]{G1982a} proved
inequality \eqref{E2.8} for real-valued functions
by a different method. Other special cases of
inequality \eqref{E2.8}
established earlier include
$K=Q^m,\,V=\sa Q^m,\,M(K,V)=\sa$ (see Bernstein \cite{B1948});
$K=\BB^m,\,V=\sa\BB^m,\,M(K,V)=\sa$
(see Nikolskii \cite[Sect. 3.2.6]{N1969}); and
$K=\BB^m,\,M(K,V)=d(V)/2$
by \eqref{E1.2a} (see the author \cite[Theorem 4]{G1979}).
\end{remark}

\begin{remark}\label{R2.4}
Since $\TT_{\sa V}\subset B_{\sa V},\, \sa>0$, inequalities
\eqref{E2.4}, \eqref{E2.8},
and \eqref{E2.9} hold for
trigonometric polynomials from $\TT_{\sa V}$.
Moreover, these inequalities are asymptotically sharp
in $\TT_{\sa V}$ as $\sa\to\iy$.
For example, the function
$f_0(x):=C_1 e^{i(k_0,x)}+C_2 e^{-i(k_0,x)}$,
where $k_0\in \sa V\cap\Z^m$ is a closest point to
$\sa \mathbf{a}\in\sa V$, belongs to $\TT_{\sa V}$, and
\ba
\sup_{x\in\R^m}
\left\|\nabla f_0(x)\right\|_K
\ge (\sa M(K,V)-C)\left\|f_0\right\|_{C(\R^m)},
\ea
where $C>0$ is a constant independent of $\sa$.
\end{remark}

\begin{remark}\label{R2.4a}
Versions of inequalities \eqref{E2.4} and \eqref{E2.8}
for continuous $n$-homogeneous polynomials
on a real or complex Hilbert space were obtained
by Anagnostopoulos, Sarantopoulos, and Tonge
\cite[Theorem 2.2]{AST2012}.
\end{remark}

\begin{example}\label{Ex2.5}
If $K=V_\mu$ and
$V=V_\la, \,1\le \mu,\,\la\le\iy$, then
\beq\label{E2.12}
M\left(V_\mu,V_\la\right)=\left\{\begin{array}{ll}
m^{1/\mu-1/\la},&1\le \mu<\la\le\iy,\\
1,&1\le \la\le\mu\le\iy.\end{array}\right.
\eeq
Therefore, by Theorem \ref{T2.1} and Corollary \ref{C2.2},
sharp inequalities \eqref{E2.4}, \eqref{E2.8},
and \eqref{E2.9} hold with
$M\left(V_\mu,V_\la\right)$ defined by \eqref{E2.12}.
In addition,
\ba
\mathbf{a}=\mathbf{a}\left(V_\mu,V_\la\right)
=\left\{\begin{array}{ll}
\left(\pm m^{-1/\la},\ldots,\pm m^{-1/\la}\right),
&1\le \mu<\la\le\iy,\\
(0,\ldots,\pm 1,\ldots,0),
&1\le \la\le\mu\le\iy,\end{array}\right.
\ea
\ba
\mathbf{b}
=\mathbf{b}\left(\left(V_\la\right)^*,
\left(V_\mu\right)^*\right)
&=&\left\{\begin{array}{ll}
\left(\pm m^{1/\mu-1},\ldots,\pm m^{1/\mu-1}\right),
&1\le \mu<\la\le\iy,\\
(0,\ldots,\pm 1,\ldots,0),
&1\le \la\le\mu\le\iy,\end{array}\right.\\
&=&\pm\left(M\left(V_\mu,V_\la\right)\vert
\mathbf{a}\vert^{-2}\right)
\mathbf{a}.
\ea
\end{example}

\section{Extremal Functions}\label{S3}
\setcounter{equation}{0}
\noindent
Throughout the section points
$\mathbf{a}\in \A(K,V)$ and
$\mathbf{b}\in \B\left(V^*,K^*\right)$
 are defined by
\eqref{E1.1a}, \eqref{E1.1b}, and \eqref{E1.2}.

Theorem \ref{T2.1} and Corollary \ref{C2.2} show that
equalities hold in \eqref{E2.4} (or \eqref{E2.8})
 and \eqref{E2.9} for functions \eqref{E2.5} and \eqref{E2.9a},
 respectively.
 In particular, functions \eqref{E2.9a} are \emph{extremal}
in inequalities \eqref{E2.8} and \eqref{E2.9}, i.e.,
there exists $x=x_0\in\R^m$ such that
equality holds in \eqref{E2.8} and \eqref{E2.9}
for functions \eqref{E2.9a}.

However, unlike the univariate case, a real-valued
 extremal function
in \eqref{E2.8} and \eqref{E2.9}
for $m\ge 2$ does not necessarily coincide
with \eqref{E2.9a}. For example,
the function $f_0(x):=R\cos \sa\vert x\vert,\,R\in\R,\,\sa>0$,
 belongs to
$B_{\sa \BB^m}$, and it is easy to verify that $f_0$
is an extremal function in \eqref{E2.8} and \eqref{E2.9}
for $V=\sa \BB^m$,
while $f_0$ cannot be represented in the form of \eqref{E2.9a}.

Nevertheless, the following theorem shows that there
are certain relations between functions \eqref{E2.9a} and
real-valued extremal functions
in \eqref{E2.8} and \eqref{E2.9}.

\begin{theorem}\label{T3.1n}
Let $f_0\in B_V\cap C(\R^m)$ be a real-valued extremal function
in \eqref{E2.8} or \eqref{E2.9}.
Then  there exists a function
$g(x):=R_1 \cos \,(\mathbf{a},x)+R_2 \sin \,(\mathbf{a},x)$,
where
$\mathbf{a}\in\A,\,R_1\in\R^1$, and $R_2\in\R^1$,
and there exists a straight line $x=x_0+\tau y_0$
in $\R^m$, where
$x_0\in\R^m,\,y_0\in\B$ and $\tau\in\R^1$,  such that
for all
$\tau\in\R^1$ the following equalities hold:
\bna
&&f_0(x_0+\tau y_0)=g(x_0+\tau y_0),\label{E3.1n}\\
&&\nabla f_0(x_0+\tau y_0)
=\nabla g(x_0+\tau y_0).\label{E3.2n}
\ena
In addition, $f_0(x_0)=0$.
\end{theorem}
\proof
\textbf{Step 1.}
We first note that if $f_0(\cdot)\in B_V\cap C(\R^m)$
is a real-valued extremal function in \eqref{E2.8},
then $f_0\left(\cdot+u\right)$ is an
extremal function in \eqref{E2.9}
for any fixed $u\in\R^m$.
So without loss of generality
we can assume that
$\left\|f_0\right\|_{C(\R^m)} =1$ and
\beq\label{E3.3n}
\left\|\nabla f_0(0)\right\|_K
=M(K,V)\sqrt{1-f_0^2(0)},
\eeq
i.e., it suffices to prove \eqref{E3.1n}
and \eqref{E3.2n} for $x_0=0$.
\vspace{.12in}\\
\textbf{Step 2.}
Next, we prove \eqref{E3.1n}. Setting
\beq\label{E3.4n}
\mathbf{a}:=\left(1/\sqrt{1-f_0^2(0)}\right)\nabla f_0(0),
\eeq
we see from \eqref{E2.9} for $K=V$ that $\mathbf{a}\in V$.
In addition, $\|\mathbf{a}\|_K=M(K,V)$ by \eqref{E3.3n}.
Therefore, $\mathbf{a}\in\A$, and
by Proposition \ref{P1.1},
there exists $\mathbf{b}=y_0\in\B$ such that
$M(K,V)=M(V^*,K^*)=\vert(\mathbf{a},\mathbf{b})\vert$.

Then setting $\vphi(\tau):=f_0\left(\tau y_0\right),
\,\tau\in\R^1$,
we see from Lemma \ref{L2.1b}
that $\vphi$ is a univariate entire function
of exponential type $\left\|y_0\right\|_{V^*}=M(K,V)$
and $\vphi\in C(\R^1)$.
Thus
\beq\label{E3.5n}
\frac
{\left\|\nabla f_0(0)\right\|_K}
{\sqrt{1-f_0^2(0)}}
=\frac
{\left\vert\left(\nabla f_0(0),y_0\right)\right\vert}
{\sqrt{1-f_0^2(0)}}
=\frac
{\left\vert\vphi^\prime(0)\right\vert}
{\sqrt{1-\vphi^2(0)}}
=M(K,V)
=M(K,V)\|\vphi\|_{C(\R^1)},
\eeq
where the last equality in \eqref{E3.5n} follows from \eqref{E2.3}.
Then
there exist $R_1\in\R^1$ and $R_2\in\R^1$ such that
$\vphi(\tau)=R_1 \cos\,[M(K,V)\,\tau]
+R_2 \sin\,[M(K,V)\,\tau]$.
Taking into account the relations
$\left\|\vphi^\prime\right\|_{C(\R^1)}
=\left\vert\vphi^\prime(0)\right\vert$
and $\|\vphi\|_{C(\R^1)}=1$,
we arrive at the equality
$\vphi(\tau)=R_2\sin\,[M(K,V)\,\tau]$,
where
$R_2:=\mbox{sgn}
\left(\nabla f_0(0),y_0\right)$.
Therefore, $f_0(0)=\vphi(0)=0$
and $g(x):=R_2\sin\, (\mathbf{a},x)$, where
$\mathbf{a}:=\nabla f_0(0)$
by \eqref{E3.4n}.
This proves \eqref{E3.1n}.\vspace{.12in}\\
\textbf{Step 3.}
Furthermore, we prove \eqref{E3.2n}.
Without loss of generality we can assume that
$\left(\nabla f_0(0),y_0\right)>0$, i.e.,
 $g(x)=\sin\, (\mathbf{a},x)$. Next,
 let us define the function
 $h:=f_0-g$
 and the straight line
 $\mathbf{L}:=\{x\in\R^m:x=\tau y_0,\tau\in\R^1\}$.
 It is clear that $h\in B_V\cap C(\R^m)$.
 Then
\beq\label{E3.6n}
h(x)=0,\qquad x\in\mathbf{L},
\eeq
by \eqref{E3.1n} and
\beq\label{E3.7n}
\nabla h(0)=0,
\eeq
since by the definition of $\mathbf{a}$ and $g$,
$\nabla f_0(0)=\mathbf{a}=\nabla g(0)$. Next, let
\ba
\tau_k:=\frac{(2k+1)\pi}{2M(K,V)},\qquad U_k:=\tau_k y_0,
\qquad k=0,\pm 1,\ldots,
\ea
be points on $\R^1$ and $\mathbf{L}$, respectively.
In addition, let $H_k$ be the hyperplane passing
through $U_k,\,k=0,\pm 1,\ldots,$
that is orthogonal to $\textbf{a}$. Note that
$U_k$ is the only element of $H_k\cap\mathbf{L},
\,k=0,\pm 1,\ldots,$
since $\left(\mathbf{a},y_0\right)>0$.

Then $g$ is constant on $H_k$ and
\ba
g(x)=g\left(U_k\right)=\sin \tau_k=(-1)^k,
\qquad x\in H_k,\quad k=0,\pm 1,\ldots.
\ea
Therefore,
\ba
(-1)^kf_0(x)\le \left\|f_0\right\|_{C(\R^m)}
=(-1)^kg(x),\qquad
x\in H_k,\quad k=0,\pm 1,\ldots,
\ea
which implies
\beq\label{E3.8n}
(-1)^kh(x)\le 0
,\qquad
x\in H_k,\quad k=0,\pm 1,\ldots.
\eeq
Since by \eqref{E3.6n},
$h\left(U_k\right)=0,\,k=0,\pm 1,\ldots,$
we see from \eqref{E3.8n}
that $h$ has a relative maximum (minimum)
on $H_k$ at $U_k$ for even (odd) numbers
$k=0,\pm 1,\ldots$. Choosing now $m-1$
linearly independent vectors
$\left\{y^{(d,k)}\right\}_{d=1}^{m-1}$
in $H_k$, we obtain by a necessary condition
for a relative extremum that
\beq\label{E3.9n}
\left(\nabla h(U_k),y^{(d,k)}\right)=0,\qquad
1\le d\le m-1,\quad k=0,\pm 1,\ldots.
\eeq
In addition, using \eqref{E3.6n}, we have that
\beq\label{E3.10n}
\left(\nabla h(U_k),y_0\right)=0,
\qquad k=0,\pm 1,\ldots.
\eeq
Then the system of vectors
$\left\{y_0, y^{(1,k)},\ldots,y^{(m-1,k)}\right\}$
is linearly independent (we recall that
$H_k\cap\mathbf{L}=\left\{U_k\right\},
k=0,\pm 1,\ldots$), and
combining \eqref{E3.9n} with \eqref{E3.10n},
we arrive at the equality
\beq\label{E3.11n}
\nabla h(U_k)=0,\qquad k=0,\pm 1,\ldots.
\eeq
Next, $h_j(x):=\partial h(x)/\partial x_j\in C(\R^m)$
by \eqref{E2.8} and, in addition,
$h_j\in B_V,\,1\le j\le m$
(see \cite[Sect. 3.1]{N1969} and \cite[Lemma 2.1 (d)]
{G2020}).
Then Step 2 of the proof of Theorem \ref{T2.1}
shows that
 the restriction $\vphi_j$ of $h_j$ to
$\mathbf{L}$ belongs to
$B_{M(K,V)}\cap C(\R^1),\,1\le j\le m$.
We also take account of the equations
$\vphi_j(\tau_k)=0,\,k=0,\pm 1,\ldots,\,1\le j\le m$,
that follow from  \eqref{E3.11n}.

Therefore, by a well-known result (see, e.g.,
\cite[Sect. 4.3.1]{T1963}),
$\vphi_j(\tau)=C_j \cos\,[M(K,V)\,\tau]$,
and taking account of $\vphi_j(0)=0$ that follows from
\eqref{E3.7n}, we obtain that for all $\tau\in\R^1,\,
\vphi_j(\tau)=0,\,1\le j\le m$.
Thus \eqref{E3.2n} is established.
\hfill $\Box$

\begin{remark}\label{R3.2n}
Theorem \ref{T3.1n} is an analogue of properties
 of extremal polynomials in Markov--type inequalities
 proved by Kro\'{o}
\cite[Corollary 1]{K2001} and R\'{e}v\'{e}sz
\cite[Theorem 1]{R2001},
and certain ideas from \cite{K2001, R2001} are used
in the proof of Theorem \ref{T3.1n}.
We do not know as to whether the corresponding analogues of
the criteria for the uniqueness of extremal polynomials
 discussed in
\cite[Theorem 2]{K2001} and
\cite[Theorem 2]{R2001} are valid for extremal functions
in Bernstein--type inequalities.
\end{remark}

\section{Markov-type Inequalities}\label{S4}
\setcounter{equation}{0}
\noindent
Throughout the section  a point $\mathbf{a}^*
  =\left(\mathbf{a}^*_1,\ldots,\mathbf{a}^*_m\right)
  \in \partial\left(V^*\right)$
   satisfies the equality
$\left\| \mathbf{a}^*\right\|_K=M(K,V^*)$,
i.e., $\mathbf{a^*}\in \A(K,V^*)$
(see \eqref{E1.1} and \eqref{E1.1a}),
and $I^m:=\left\{u\in\R^m:
 u_j\ge 0,\,1\le j\le m\right\} $
 is the first orthant of $\R^m$.

Markov-type inequalities on convex bodies with sharp constants
have been a hot topic in approximation theory since the early 1990s
(see surveys \cite{K2005} and \cite{CS2019}).
\vspace{.12in}\\
\textbf{Symmetric Bodies.}
The following multivariate version of the Markov inequality
for polynomials with real coefficients
was  proved by Sarantopoulos \cite[Theorem 2]{S1991}
who elegantly applied the inscribed ellipse method
(IEM) to this problem.

\begin{theorem}\label{T3.1a}
For $P\in\PP_{n,m}$,
\beq\label{E3.1}
\left\vert\left(\nabla P(x),y\right)
\right\vert\le \| y\|_V\,n^2 \,\|P\|_{C(V)},
\qquad x\in V,\quad y\in\R^m\setminus \{0\}.
\eeq
\end{theorem}
\noindent
Independently and by the pluripotential theory approach (PTA), Baran
\cite[Theorem 2 (c)]{B1994} established \eqref{E3.1} and extended
this inequality to any $P\in\PP_{n,m}$ with complex coefficients.
It turned out later that IEM and PTA are equivalent
(see Burns, Levenberg, Ma'u, and R\'{e}v\'{e}sz
 \cite[Corollary 4.3]{BLMR2010}).
In case of $V=\BB^m$, \eqref{E3.1} was earlier obtained by
Kellogg \cite[Theorem VI]{K1928}.

\begin{remark}\label{R3.1b}
For the reader's convenience, we present here a shorter
 and more straightforward
 proof of \eqref{E3.1}, compared with \cite{B1994},
  for complex-valued polynomials
  (cf. \cite[Theorem 1.1]{GT2017}).
  Let $y\in\R^m\setminus\{0\}$  and let
  $P=P_1+iP_2\in \PP_{n,m}$ be a nonconstant
polynomial, where $P_j\in \PP_{n,m},\,j=1,\,2,$ are
polynomials with real coefficients.
Then there exists $\xi\in V$ such that
$\max_{x\in V}
\left\vert\left(\nabla P(x),y\right)\right\vert
=\left\vert\left(\nabla P(\xi),y\right)\right\vert$.
Let us define $\g\in  [0,2\pi)$ by the equality
$e^{i\g}=\left(\nabla P(\xi),y\right)
/\left\vert \left(\nabla P(\xi),y\right)\right\vert$.
Then the polynomial
$D:=\cos \g\,P_1+\sin \g\,P_2\in\PP_{n,m}$
has real coefficients
and satisfies the relations:
$\left\vert D(x)\right\vert\le \left\vert P(x)\right\vert$
for $x\in V$ and
$\left\vert \left(\nabla D(\xi),y\right)\right\vert
= \left\vert \left(\nabla P(\xi),y\right)\right\vert$.
Therefore, for any $x\in V$ and $y\in\R^m\setminus\{0\}$,
\ba
 \frac{\left\vert \left(\nabla P(x),y
 \right)\right\vert}{\|P\|_{C(V)}}
 \le \frac{\left\vert \left(\nabla P(\xi),y
 \right)\right\vert}{\|P\|_{C(V)}}
 \le \frac{\left\vert \left(\nabla D(\xi),y
 \right)\right\vert}{\|D\|_{C(V)}}
 \le  \| y\|_V\,n^2,
 \ea
 by \cite[Theorem 2]{S1991}. Thus \eqref{E3.1} is established.
\end{remark}

Inequality \eqref{E3.1} along with \eqref{E1.1} and \eqref{E1.2a}
immediately imply the following relations:
\beq\label{E3.2}
\left\vert\nabla P(x)\right\vert
=
\left\|\nabla P(x)\right\|_{\BB^m}
\le  M(\BB^m,V)\, n^2\,\|P\|_{C(V)}
=\frac{2n^2}{w(V)}\|P\|_{C(V)},\qquad x\in V,
\eeq
(see, e.g., \cite[Eq. (7)]{KR1999}). The polynomial
$P_0(x):=CT_n((\mathbf{a}^*,x)),\,C\in\CC^1,$
is an \textit{extremal} polynomial in \eqref{E3.2},
i.e., equality holds in \eqref{E3.2} for $P=P_0$
and for $x=x_0\in V$, satisfying the
equality $(\mathbf{a}^*,x_0)=\pm 1$.

However, unlike the univariate case, a real-valued
 extremal polynomial in \eqref{E3.2}
for $m\ge 2$ does not necessarily coincide
with $P_0$. For example,
the polynomial $P(x):=R\,T_n(\vert x\vert),\,R\in\R^1,$
 for an even $n$
is an extremal polynomial in \eqref{E3.2} for $V=\BB^m$
(not $T_{n/2}\left(\vert x\vert^2\right)$ as stated in
 \cite[p. 466]{R2001}),
while $P$ cannot be represented in the form of $P_0$.

Kro\'{o}
\cite{K2001} and R\'{e}v\'{e}sz
\cite{R2001} studied properties of the extremal
polynomials in \eqref{E3.2}. In particular, they found
certain relations (like \eqref{E3.1n} and \eqref{E3.2n})
between polynomials $P_0$ for $V=\BB^m$ and
real-valued extremal functions in \eqref{E3.2}
(see \cite[Corollary 1]{K2001} and \cite[Theorem 1]{R2001}).
In addition, they provided criteria for the uniqueness of
 the extremal polynomials in \eqref{E3.2}
(see \cite[Theorem 2]{K2001} for $m=2$ and
\cite[Theorem 2]{R2001} for $m>2$).

The following result is a polynomial version of \eqref{E2.8}
 and a more general version of \eqref{E3.2}.

 \begin{corollary}\label{C3.1}
 For every $P\in \PP_{n,m}$ and $x\in V$,
the following inequality holds:
\beq\label{E3.3}
\left\|\nabla P(x)\right\|_K
\le M(K,V^*)\,n^2\,\|P\|_{C(V)}.
\eeq
 Equality holds in \eqref{E3.3} for
 $P(x)=P_0(x):=CT_n((\mathbf{a}^*,x)),\,C\in\CC^1,$
 and $x=x_0\in V\cap H_\pm(\mathbf{a}^*)$,
 where $H_\pm(\mathbf{a}^*)
 :=\{x\in\R^m:(\mathbf{a}^*,x)=\pm 1\}$ are
 parallel supporting hyperplanes of $V$.
 Hence inequality \eqref{E3.3} is sharp.
 \end{corollary}
 \proof
 It follows from Theorem \ref{T3.1a}
 that for every $y\in K^*$,
 \beq\label{E3.4}
\left\vert\left(\nabla P(x),y\right)
\right\vert\le  \| y\|_V/
\| y\|_{K^*}\,n^2\,\|P\|_{C(V)}
\le M(V,K^*)\,n^2\,\|P\|_{C(V)}
=M(K,V^*)\,n^2\,\|P\|_{C(V)},
\eeq
by \eqref{E1.2}. Then \eqref{E3.3} follows from \eqref{E3.4}.
Finally,
$
\left\|\nabla P_0(x_0)\right\|_K
=\left\| \mathbf{a}^*\right\|_{K}\,n^2\,\|P_0\|_{C(V)}.
$
\hfill $\Box$
\begin{remark}\label{R3.1c}
Note that for $K=\BB^m$ \eqref{E3.3} is reduced to \eqref{E3.2}.
In addition,
we conjecture that using the techniques from
\cite{K2001} and \cite{R2001}, it is possible
to extend the properties of the extremal polynomials in
\eqref{E3.2} to the extremal polynomials in
\eqref{E3.3}.
\end{remark}
\noindent
\textbf{Non-symmetric Bodies.}
Bia{\l}as-Cie\.{z} and Goetgheluck \cite[Sect. 3]{BG1995}
found out that \eqref{E3.2} fails for certain
non-symmetric convex bodies, in particular, for
the simplex
\ba
\Delta_0:=\left\{u\in I^m:
0\le \sum_{j=1}^m u_j\le 1\right\}.
\ea
Certain estimates of the sharp constant in the
Markov-type inequality on triangles in $\R^2$
 were obtained by Kro\'{o} and R\'{e}v\'{e}sz
\cite[Theorem 4]{KR1999}.
Other Markov-type inequalities on convex bodies
were obtained by Wilhelmsen \cite{W1974} and
Ditzian \cite{D1992}.
Skalyga \cite[Theorems 1 and 2]{S1997} proved
 that for every convex body
 $\C\subset \R^m$ and $Q\in\PP_{n,m}$,
\ba
\left\|\,\left\vert\nabla Q\right\vert\,\right\|_{C(\C)}
\le \frac{2}{w(\C)}n\,\cot\left(\frac{\pi}{4n}\right)
\left\|Q\right\|_{C(\C)} ,
\ea
and this inequality cannot be improved on the class of
all convex bodies in $\R^m$.
The same result was
independently announced by
Subbotin and Vasiliev \cite{VS1998} as well.

However, no sharp constant in the Markov-type inequality
 on the given
non-symmetric convex body $\C$ is known. Below we find
sharp constants for certain bodies
$\C$ if the gradient $\nabla f(u)$
is replaced by the "weighted" gradient vector
\beq\label{E3.5}
 \nabla_u f(u)
 :=\left(\sqrt{u_1}\frac{\partial f(u)}{\partial u_1},\cdots,
 \sqrt{u_m}\frac{\partial f(u)}{\partial u_m}\right),
 \qquad u\in I^m.
 \eeq
 We first need certain conditions on $K$ and $V$.

 \begin{definition}\label{D3.2}
 We say that a body $V\subset \R^m$ satisfies
 the $S$-condition if
 $V$ is symmetric about all coordinate hyperplanes, that is,
 for every $x\in V$ the vectors $(\pm\vert x_1\vert,\ldots,
 \pm\vert x_m\vert)$ belong to $V$.
 \end{definition}
 \noindent
 The $S$-condition is equivalent to the $\Pi$-condition
 introduced in \cite[Definition 1.1]{G2021}.

 \begin{definition}\label{D3.3}
 We say that a pair $(K,V),\,K\subset \R^m,\,V\subset \R^m$,
  satisfies the $\mathbf{a}^*$-condition
 if there exists a point
 $\mathbf{a}^*\in \A(K,V^*)$ that
 belongs to a coordinate axis.
 \end{definition}
 \noindent
 Next, given $V\subset \R^m$, we define the body
 $\C=\C(V)\subset I^m$
 (not necessarily convex) by the formula
 \beq\label{E3.6}
 \C:=\left\{u\in I^m:\left\|\left(\sqrt{u_1},
 \ldots,\sqrt{u_m}\right)\right\|_V\le 1\right\}.
 \eeq
Then the following result is valid.

 \begin{theorem}\label{T3.4}
 Let $V$ satisfy the $S$-condition and let
 $\nabla_u$ and $\C$ be defined by \eqref{E3.5} and
  \eqref{E3.6}, respectively.
 Then for every $Q\in \PP_{n,m}$ and $u\in \C$,
the following inequality holds:
\beq\label{E3.7}
\left\|\nabla_u Q(u)\right\|_K
\le 2M(K,V^*)\,n^2\,\|Q\|_{C(\C)}.
\eeq
 If $(K,V)$ satisfies
 the $\mathbf{a}^*$-condition, then
 there exists $k,\,1\le k\le m$, such that
 equality holds in \eqref{E3.7} for
 $Q(u)=Q_0(u):=CT_{2n}\left(\mathbf{a}^*_k\,
 \sqrt{u_k}\right),\,C\in\CC^1,$
 and $u=u_0:=\left(0,\ldots,0,
 \pm \left(\mathbf{a}^*_k\right)^{-2},0,\ldots,0\right)\in \C$.
 \end{theorem}
 \proof
 First note that the following statement
 is an immediate consequence of
 Definition \ref{D3.2} and the definition of $\C$:
 $x\in V$ if and only if
 $u=\left(x_1^2,\ldots,x_m^2\right)\in\C$.
Next, by Corollary \ref{C3.1}, for every
 $P(x)
 =\sum_{\vert\al\vert\le n}c_\al x^{2\al}\in\PP_{2n,m}$
 and $x\in V$, the following inequality holds:
\beq\label{E3.8}
\left\|\nabla P(x)\right\|_K
\le 4 M(K,V^*)\,n^2\,\|P\|_{C(V)}.
\eeq
Therefore, it follows from \eqref{E3.8}
and the aforementioned statement
that for every polynomial
$Q(u)
=\sum_{\vert\al\vert\le n}c_\al u^{\al}\in\PP_{n,m}$
and $u\in \C$,
\ba
\left\|\nabla_u Q(u)\right\|_K
=(1/2)\left\|\nabla P(x)\right\|_K
\le 2 M(K,V^*)\,n^2\,\|P\|_{C(V)}
=2 M(K,V^*)\,n^2\,\|Q\|_{C(\C)}.
\ea
Thus \eqref{E3.7} follows. If $(K,V)$ satisfies the
$\mathbf{a}^*$-condition, then there exists
$k,\,1\le k\le m$, such that
$\mathbf{a}^*
=\left(0,\ldots,0,\mathbf{a}^*_k,0,\ldots,0\right)$.
Then by Corollary \ref{C3.1},
equality holds in \eqref{E3.8} for
$P(x)=P_0(x):=CT_{2n}\left(\mathbf{a}^*_k x_k\right)
\in\PP_{2n,m},\,C\in\CC^1,$
 and $x=\left(0,\ldots,0,\pm 1/\mathbf{a}^*_k,0,\ldots,0\right)
 \in V\cap H_\pm(\mathbf{a}^*)$,
 where $H_\pm(\mathbf{a}^*)
 :=\left\{x\in\R^m:x_k=\pm 1/\mathbf{a}^*_k\right\}$ are
 parallel supporting hyperplanes of $V$.
 Therefore, equality holds in \eqref{E3.7}
 for $Q(u)=Q_0(u):=CT_{2n}\left(\mathbf{a}^*_k \sqrt{u_k}\right)
 \in\PP_{n,m}$
 and $u=u_0:=\left(0,\ldots,0,
 \pm \left(\mathbf{a}^*_k\right)^{-2},0,\ldots,0\right)\in \C$.
\hfill $\Box$

\begin{example}\label{Ex3.5}
If $K=V_\mu$ and
$V=V_\la, \,1\le \mu,\,\la\le\iy$, then
\beq\label{E3.10}
M\left(V_\mu,\left(V_\la\right)^*\right)=\left\{\begin{array}{ll}
m^{1/\mu+1/\la-1},&1/\mu+1/\la>1,\\
1,&1/\mu+1/\la\le 1.\end{array}\right.
\eeq
In addition,
\ba
\mathbf{a}^*=\mathbf{a}^*\left(V_\mu,\left(V_\la\right)^*\right)
=\left\{\begin{array}{ll}
\left(\pm m^{1/\la-1},\ldots,\pm m^{1/\la-1}\right),
&1/\mu+1/\la>1,\\
(0,\ldots,\pm 1,\ldots,0),&1/\mu+1/\la\le 1.\end{array}\right.
\ea
Furthermore, $V_\mu$ satisfies the $S$-condition
for $1\le \mu\le\iy$,
and $\left(V_\mu,V_\la\right)$ satisfies
the $\mathbf{a}^*$-condition for $1/\mu+1/\la\le 1$.
Then by Theorem \ref{T3.4},
sharp inequality \eqref{E3.7} holds for
$K=V_\mu,\,\C=\C_\la:=V_{\la/2}\cap I^m$,
and $1/\mu+1/\la\le 1$
with $M\left(V_\mu,\left(V_\la\right)^*\right)$
defined by \eqref{E3.10}.
 Note that $\C_\la$
is not convex for $1\le \la<2$.
 In particular, for the simplex
$\Delta_0=\C_2$ and $\mu=2$ the sharp inequality
\ba
\left\vert\nabla_u Q(u)\right\vert
\le 2 n^2\,\|Q\|_{C(\Delta_0)},
\qquad Q\in\PP_{n,m},\quad u\in\Delta_0,
\ea
is valid.
\end{example}
\noindent
\textbf{Acknowledgement}
We are grateful to
Leokadia Bia{\l}as-Cie\.{z} for the
provision of references.

\end{document}